\documentclass{amsart}
\usepackage[latin1]{inputenc}
\usepackage[T1]{fontenc}
\usepackage{ae,aecompl}

\usepackage{subfigure}

\usepackage{amssymb}
\usepackage{amsmath}
\usepackage{units}

\usepackage[english]{babel}

\usepackage[normalem]{ulem}

\usepackage{tikz,ifthen}
\usetikzlibrary{matrix,shapes,arrows,chains}
\usetikzlibrary{calc}
\usetikzlibrary{patterns}
\usetikzlibrary{positioning} 

\definecolor{darkgreen}{rgb}{0,0.7,0}
\definecolor{darkred}{rgb}{0.7,0,0}
\definecolor{darkblue}{rgb}{0,0,0.7}
\usepackage{pdflscape}

\usepackage{enumerate}
\usepackage{verbatim}
\usepackage{hyperref}

\usepackage{graphicx}
\usepackage{rotating}
\usepackage{algorithmic}
\usepackage{algorithm}
\usepackage{xspace}

\usepackage{ifpdf}
\ifpdf
\fi

\newtheorem{theorem}{Theorem}[section]

\newtheorem{proposition}[theorem]{Proposition}
\newtheorem{corollary}[theorem]{Corollary}
\newtheorem{definition}[theorem]{Definition}

\newtheorem{problem}[theorem]{Problem}

\newcommand{\x}{\mathbf{x}}
\newcommand{\sg}[1][n]{\mathfrak{S}_{#1}}

\newcommand{\QQ}{\mathbb{Q}}

\newcommand{\KK}{\mathbb{K}}

\newcommand{\Sage}{\texttt{Sage}\xspace}
\newcommand{\sagecombinat}{\texttt{Sage-Combinat}\xspace}
\newcommand{\gap}{\texttt{GAP}\xspace}

\newcommand{\ind}{\hspace{4ex}}

\definecolor{sagecolor}{rgb}{.78, .36, .04}
\newcommand{\sageex}[1]{\vspace{0.5em} {\small \tt #1} \vspace{0.5em}}
\def\sageret#1{\begin{center} #1 \end{center}}

\def\sagepromt{{\normalsize \color{sagecolor} sage: }}
\def\sagedots{{\normalsize \color{sagecolor}....: }}

\usepackage{ifthen}
\newboolean{draft}
\setboolean{draft}{true}
\ifdraft
\newcommand{\TODO}[2][To do: ]{\textcolor{red}{\textbf{#1#2}}}
\newcommand{\INFO}[2][Info: ]{\textcolor{red}{\textbf{#1#2}}}
\else
\newcommand{\TODO}[2][]{}
\newcommand{\INFO}[2][]{}
\fi

\title[Effective invariant theory using representation
  theory]{Effective invariant theory of permutation groups using
  representation theory}

\author{Nicolas Borie}

\address{Univ. Paris Est
  Marne-La-Vall\'ee, Laboratoire d'Informatique Gaspard Monge, Cit\'e
  Descartes, B\^at Copernic -- 5, bd Descartes Champs sur Marne 77454
  Marne-la-Vall\'ee Cedex 2, France}

\begin{document}

% \tableofcontents

\maketitle

\begin{abstract}
  Using the theory of representations of the symmetric
  group, we propose an algorithm to compute the invariant ring of a
  permutation group. Our approach have the goal to reduce the amount
  of linear algebra computations and exploit a thinner combinatorial
  description of the invariant ring.
\end{abstract}


\begin{center}
\textbf{This is a drafty old version : full corrected text are
  available at : http://www.springer.com/}
\end{center}

%%%%%%%%%%%%%%%%%%%%%%%%%%%%%%%%%%%%%%%%%%%%%%%%%%%%%%%%%%%%%%%%%%%%%%%%%%%%%%%%
\section{Introduction}
%%%%%%%%%%%%%%%%%%%%%%%%%%%%%%%%%%%%%%%%%%%%%%%%%%%%%%%%%%%%%%%%%%%%%%%%%%%%%%%%

Invariant theory has been a rich and central area of algebra ever
since the eighteenth theory, with practical
applications~\cite[\S~5]{Kemper_Derksen.CIT.2002} in the resolution of
polynomial systems with symmetries (see
e.g. \cite{Colin.1997.SolvingSymmetries},
\cite{Gatermann.1990.Symmetry}, \cite[\S~2.6]{Sturmfels.AIT},
\cite{Faugere_Rahmany.2009.SAGBIGroebner}), in effective Galois theory
(see e.g.~\cite{Colin.TIE}, \cite{Abdeljaouad.TIATG},
\cite{Geissler_Kluners.2000.GaloisGroupComputations}), or in discrete
mathematics (see e.g.~\cite{Thiery.AIG.2000,Pouzet_Thiery.IAGR.2001}
for the original motivation of the second author). The literature
contains deep and explicit results for special classes of groups, like
complex reflection groups or the classical reductive groups, as well
as general results applicable to any group. Given the level of
generality, one cannot hope for such results to be simultaneously
explicit and tight in general. Thus the subject was effective early
on: given a group, one wants to \emph{calculate} the properties of its
invariant ring. Under the impulsion of modern computer algebra,
computational methods, and their implementations, have largely
expanded in the last twenty
years~\cite{Kemper.Invar,Sturmfels.AIT,Thiery.CMGS.2001,Kemper_Derksen.CIT.2002,King.2007.secondary,King.2007.minimal}.
However much progress is still needed to go beyond toy examples and
enlarge the spectrum of applications.

Classical approaches solving the problem of computing invariant ring
use elimination techniques in vector spaces of too high
dimensions. Gröbner basis become impracticable when the number of
variables goes up (around 10 for modern computers). The evaluation
approach proposed by the author in~\cite{Borie.2011.Thesis} required a
permutation group whose index in the symmetric group is relatively
controlled (around 1000 for modern computers). Each approaches
localize the algebra reduction in vector spaces still of too large
dimensions. Gröbner basis approaches works with monomials of degree
$d$ over $n$ variables and linear reduction over a space spanned by
these monomials is costly. The evaluation approach, proposed by the
author in his thesis, does the linear algebra in a free module spanned
by the cosets of the symmetric group by a permutation group (i.e. it
have the index for dimension). In both case, as linear reduction
globally cost the cube of the dimension of the space, one cannot hope
to go so much further with classical approaches even with the progress
of computer.

We propose in this article an approach following the idea that adding
more combinatorics in invariant theory help to produce more efficient
algorithms whose outputs could perhaps reveal some combinatorics also
; the long time goal being having a combinatorial description of
invariant ring (generators or couple of primary secondary invariants
families). Since Hilbert, This problem have been solve only in very
restrictive and special cases (for example, \cite{Garsia1984107} give
secondary invariants for Young subgroups of symmetric groups). We
focus on the problem of computing secondary invariants of finite
permutation groups in the non modular case. Assuming that, we will
shows how to localize computations inside selected irreducible
representations of the symmetric group. These spaces are smaller than
the ones used in classical approaches and we can largely take
advantages of the combinatorial results coming from the theory of
representations of the symmetric group.

\section{Invariant ring and representations of the symmetric group}

For any set $S$, we will denote $|S|$ the cardinality of the set $S$.

\subsection{Invariant ring of permutation group and application to combinatorics}

Our approach start from a result in one key article of invariant
theory written by Stanley~\cite[proposition 4.9]{Stanley.1979} mixing
invariant of finite group and combinatorics. We recall this general
result:

\begin{proposition}
Let $\theta_1, \dots , \theta_n$ be an homogeneous set of parameters
for $\KK[\x]^G$, where $G$ is any finite subgroup of $GL(\KK^n)$ of
order $|G|$. Set $d_i = deg(\theta_i)$ and $t = d_1 \dots d_n / |G|$. Then
the action of $G$ on the quotient ring $S = \KK[\x]/(\theta_1, \dots ,
\theta_n)$ is isomorphic to $t$ times the regular representation of
$G$.
\end{proposition}

Applying this result to $\sg[n]$ the symmetric group of degree $n$
with $\theta_i = e_i$ (elementary symmetric polynomial)
\begin{equation}
  e_i := \sum_{P \subset \{1, \dots , n\}} ( \prod_{i \in P} x_i ),
\end{equation}
we recover the well known result that the ring $\KK[\x] / (e_1, \dots
, e_n)$ is isomorphic the regular representation $RR(\sg[n])$ of the symmetric
group (here $t = \frac{n!}{n!} = 1$). This well known quotient $\KK[\x] / (e_1, \dots , e_n)$
is called the coinvariant ring of the symmetric group in the algebraic
combinatorics world and several basis of this ring have been
explicitly built (Harmonic polynomials, Schubert Polynomials, Descents
monomials and more~\cite{MR2538310,Lascoux-Schutzenberger.1982,Garsia1984107}).

\begin{equation}
  \KK[\x] / (e_1, \dots , e_n) \sim_{\sg[n]} RR(\sg[n])
\end{equation}

Let $G$ a group of permutations, subgroup of $\sg[n]$, we know reapply
the result of Stanley to $G$ with the same homogeneous set of
parameters formed with the elementary symmetric polynomials. Then, the
ring of coinvariant of the symmetric group is now also isomorphic to
$t = n!/|G|$ time the regular representation of the group $G$.

\begin{equation}
  \KK[\x] / (e_1, \dots , e_n) \sim_{G} \bigoplus_{i=1}^{n!/|G|} RR(G)
\end{equation}

We know that for any permutation group $G$, in the non modular case,
the ring of invariant under the action of $G$ is a Cohen-Macaulay
algebra. This imply that there exist a family of generator $\eta_i$
making the ring of invariant under the action of $G$ a free module of
rank $n!/|G|$ over the ring of symmetric polynomials.

\begin{equation}
  \KK[\x]^G = \bigoplus_{i=1}^{n!/|G|} \eta_i \KK[e_1, \dots , e_n]
\end{equation}

Taking the quotient on both side by the ideal $(e_1, \dots , e_n)$ and
keeping $\eta_i$ as representative of its equivalent class in the
quotient, we have

\begin{equation}
  \KK[\x]^G / (e_1, \dots , e_n)  = \bigoplus_{i=1}^{n!/|G|} \KK \cdot \eta_i
\end{equation}

As $\KK[\x]^G$ is, by definition, the subspace of $\KK[\x]$ on which
the action of $G$ is trivial, the result of Stanley imply in
particular that the polynomials $\eta_i$ span the subspace of the
coinvariant of the symmetric group on which the action of $G$ is
trivial. A way to construct the $\eta_i$ is thus to search them as
$G$-fixed point inside the ring of coinvariant of the symmetric group
and that could be done irreducible representation by irreducible
representation.

The theory of representations of the symmetric group has been largely
studied, this bring us to formulate the following problem:

\begin{problem}\label{main-pb}
Let $n$ a positive integer and $G$ a permutation group, subgroup of
$\sg[n]$. Construct an explicit basis of the trivial representations
of $G$ appearing in each irreducible subrepresentation of $\sg[n]$
inside the quotient $\KK[\x] / (e_1, \dots , e_n)$.
\end{problem}

A first step to solve this problem constitute in having a basis of the
coinvariant for the symmetric group respecting the action of $\sg[n]$
such that it can be partitioned by irreducible representations. We
expose this basis in the next section.

\subsection{Representations of the symmetric group}

We recall in this section some results describing the irreducible
representations of the symmetric group.

For a positive integer $n$, we will call $\lambda$ partition of $n$
(denoted $\lambda \vdash n$) a non increasing sequence of integers
$\lambda = (\lambda_1, \dots , \lambda_r)$ whose entries sum up to
$n$.

\sageex{
\sagepromt P = Partitions(4); P

\sageret{Partitions of the integer 4}

\sagepromt P.list()

\sageret{[[4], [3, 1], [2, 2], [2, 1, 1], [1, 1, 1, 1]]}

\sagepromt Partitions(8).random\_element()

\sageret{[3, 2, 2, 1]}
}

For a positive integer $n$, irreducible representations of the
symmetric group $\sg[n]$ are indexed by the Partitions of $n$. Since
we have a finite group, the multiplicity of an irreducible
representation inside the regular representation equal its
dimension. This information can be collected studying \emph{standard
  tableaux}.

Let $n$ a positive integer and $\lambda = (\lambda_1, \dots ,
\lambda_r)$ a partition of $n$. A tableau of shape $\lambda$ is a
diagram of square boxes disposed in raw such that the first raw
contains $\lambda_1$ boxes, on top of its, a second raw contains
$\lambda_2$ boxes and so on. A \emph{standard tableau} of shape
$\lambda$ is a filled tableau of shape $\lambda$ with integer from $1$
to $n$ such that integers are increasing in column and raw. We will
denote $STab(\lambda)$ the set of standard tableaux of shape
$\lambda$.

We can ask \texttt{Sage} to display a such object.

\sageex{
\sagepromt T = StandardTableaux([3,2,2,1]).random\_element(); T

\sageret{[[1, 2, 3], [4, 5], [6, 8], [7]]}

\sagepromt latex(T)

{\def\lr#1{\multicolumn{1}{|@{\hspace{.6ex}}c@{\hspace{.6ex}}|}{\raisebox{-.3ex}{$#1$}}}
\raisebox{-.6ex}{$\begin{array}[t]{*{3}c}\cline{1-1}
\lr{7}\\\cline{1-2}
\lr{6}&\lr{8}\\\cline{1-2}
\lr{4}&\lr{5}\\\cline{1-3}
\lr{1}&\lr{2}&\lr{3}\\\cline{1-3}
\end{array}$}
}
}

We can also iterate or generate all tableaux of a given shape.

\sageex{
\sagepromt S = StandardTableaux([2,2,1]); S

\sageret{Standard tableaux of shape [2, 2, 1]}

\sagepromt for T in S: latex(T)

{\def\lr#1{\multicolumn{1}{|@{\hspace{.6ex}}c@{\hspace{.6ex}}|}{\raisebox{-.3ex}{$#1$}}}
\raisebox{-.6ex}{$\begin{array}[t]{*{2}c}\cline{1-1}
\lr{3}\\\cline{1-2}
\lr{2}&\lr{5}\\\cline{1-2}
\lr{1}&\lr{4}\\\cline{1-2}
\end{array}$}
}
{\def\lr#1{\multicolumn{1}{|@{\hspace{.6ex}}c@{\hspace{.6ex}}|}{\raisebox{-.3ex}{$#1$}}}
\raisebox{-.6ex}{$\begin{array}[t]{*{2}c}\cline{1-1}
\lr{4}\\\cline{1-2}
\lr{2}&\lr{5}\\\cline{1-2}
\lr{1}&\lr{3}\\\cline{1-2}
\end{array}$}
}
{\def\lr#1{\multicolumn{1}{|@{\hspace{.6ex}}c@{\hspace{.6ex}}|}{\raisebox{-.3ex}{$#1$}}}
\raisebox{-.6ex}{$\begin{array}[t]{*{2}c}\cline{1-1}
\lr{4}\\\cline{1-2}
\lr{3}&\lr{5}\\\cline{1-2}
\lr{1}&\lr{2}\\\cline{1-2}
\end{array}$}
}
{\def\lr#1{\multicolumn{1}{|@{\hspace{.6ex}}c@{\hspace{.6ex}}|}{\raisebox{-.3ex}{$#1$}}}
\raisebox{-.6ex}{$\begin{array}[t]{*{2}c}\cline{1-1}
\lr{5}\\\cline{1-2}
\lr{2}&\lr{4}\\\cline{1-2}
\lr{1}&\lr{3}\\\cline{1-2}
\end{array}$}
}
{\def\lr#1{\multicolumn{1}{|@{\hspace{.6ex}}c@{\hspace{.6ex}}|}{\raisebox{-.3ex}{$#1$}}}
\raisebox{-.6ex}{$\begin{array}[t]{*{2}c}\cline{1-1}
\lr{5}\\\cline{1-2}
\lr{3}&\lr{4}\\\cline{1-2}
\lr{1}&\lr{2}\\\cline{1-2}
\end{array}$}
}
}

The number $f^{\lambda}$ of standard tableaux of a given shape
$\lambda$ can be easily computed using the hook-length
formula~\cite[formula 2.1]{MR2538310}. As standard tableaux of shape
$\lambda$ constitute a basis indexing of the vector space associated
to the irreducible representation of the symmetric group indexed by
$\lambda$, and because this same representation must have for
multiplicity its dimension inside the regular representation of
$\sg[n]$, we have
\begin{equation}\label{sum_f_lambda}
  \sum_{\lambda \vdash n} (f^{\lambda})^2 = |\sg[n]| = n!
\end{equation}

The following computation illustrate this equality and check that
hook-length formula is well implemented in \texttt{Sage}.

\sageex{
\sagepromt def check\_multiplicities(n):

\sagedots \qquad dim = 0

\sagedots \qquad for p in Partitions(n):

\sagedots \qquad \qquad dim = dim + StandardTableaux(p).cardinality()**2

\sagedots \qquad return dim

\sagepromt for i in range(1, 9) : print i, " : ", check\_multiplicities(i)

\sageret{1  :  1 \\
2  :  2 \\
3  :  6 \\
4  :  24 \\
5  :  120 \\
6  :  720 \\
7  :  5040 \\
8  :  40320} 
}

We know describe a last useful object for our algorithmic to come
which gather information about all irreducible representations of the
symmetric group : the \emph{character table}.

We recall that character of a representation is a map which associate
the trace of matrices for each group element. Such map are constant on
conjugacy classes and conjugacy classes of the symmetric group are
indexed also by partitions. Any permutation $\tau \in \sg[n]$ have a
single disjoint cycles representation and it belong to the conjugacy
class indexed by the partition $\mu = (\mu_1, \dots , \mu_r)$ if its
disjoint cycles representation contains $r$ cycles of size
respectively $\mu_1$, $\mu_2$, $\dots$ , $\mu_r$.

The \emph{character table} of the symmetric group gather in a square
matrices the value of characters of irreducible representations on
each conjugacy classes of the symmetric group.

\sageex{
\sagepromt G = SymmetricGroup(3); G

\sageret{Symmetric group of order $3!$ as a permutation group}

\sagepromt G.character\_table()

\sageret{
\begin{displaymath}
\left(\begin{array}{rrr}
1 & -1 & 1 \\
2 & 0 & -1 \\
1 & 1 & 1
\end{array}\right)
\end{displaymath}
}

\sagepromt for C in G.conjugacy\_classes(): 

\sagedots \qquad Permutation(C.representative()).cycle\_type()

\sageret{
\texttt{[ 1, 1, 1 ]} \\
\texttt{[ 2, 1 ]} \\
\texttt{[ 3 ]}
}
}

\section{Higher Specht polynomials for the symmetric group}

Algorithmic in invariant theory must, at some point, construct
invariant polynomials. Most current approaches use the Reynolds
operator or an orbit sum over a group of a special monomial. When the
group become large, such invariant become very large even they are
stored in a sparse manner inside a computer, the number of terms can
easily fit with $n!$ when $G$ is a permutation group with a small
index in $\sg[n]$.

For our approach, with focuses on the combinatorics of the quotient
$\KK[\x] / (e_1, \dots , e_n)$, the \emph{higher specht polynomials}
will constitute the perfect family to get explicit answer to
Problem~\ref{main-pb}.

The quotient $\KK[\x] / (e_1, \dots , e_n)$ is isomorphic to the
regular representation of $\sg[n]$ in which we have several copies of
irreducible representation following their dimension. The \emph{Specht
  polynomials}, which are associated standard tableaux, allows to
construct explicit subspace of $\KK[\x]$ isomorphic to an irreducible
representation of the symmetric group. For $\lambda \vdash n$ a
partition, the span of Specht polynomials associated to standard
tableaux of shape $\lambda$ realize explicitly the irreducible
representation of $\sg[n]$ indexed by $\lambda$. Now, we will see that
the higher Specht polynomial take care of multiplicities of
irreducible representation inside the coinvariant. They are indexed
by pair of standard tableaux of the same shape and they constitute a
basis of the $\sg[n]$-module $\KK[\x] / (e_1, \dots , e_n)$. Among all
known basis of the coinvariants for the symmetric group (Harmonic,
Schubert, monomials under the staircase, descents monomials, ...), the
higher Specht polynomial constitute a basis which can be partitioned
by irreducible $\sg[n]$-module by construction.

let $\lambda \vdash n$ a partition and $S, T$ two standard tableaux of
shape $\lambda$. We define the word $w(S)$ by reading the tableau $S$
from the top to the bottom in consecutive columns, starting from the
left. The number $1$ in the word $w(S)$ as for index $0$. Now,
recursively, if the number $k$ in the word has index $p$, then $k+1$
has index $p+1$ if it lies to the left of $k$ in the word, it has
index $p$ otherwise. For example, with the two tableaux
\begin{displaymath}
S = {\def\lr#1{\multicolumn{1}{|@{\hspace{.6ex}}c@{\hspace{.6ex}}|}{\raisebox{-.3ex}{$#1$}}}
\raisebox{-.6ex}{$\begin{array}[]{*{3}c}\cline{1-2}
\lr{3}&\lr{5}\\\cline{1-3}
\lr{1}&\lr{2}&\lr{4}\\\cline{1-3}
\end{array}$}
}
\qquad
T = {\def\lr#1{\multicolumn{1}{|@{\hspace{.6ex}}c@{\hspace{.6ex}}|}{\raisebox{-.3ex}{$#1$}}}
\raisebox{-.6ex}{$\begin{array}[]{*{3}c}\cline{1-2}
\lr{2}&\lr{4}\\\cline{1-3}
\lr{1}&\lr{3}&\lr{5}\\\cline{1-3}
\end{array}$}
}
\end{displaymath}
The reading of the Tableau $S$ give $31524$, now placing step by step
the indices, we get
\begin{displaymath}
\begin{array}{cccc}
  & & 3_{\ }1_{0}5_{\ }2_{\ }4_{\ } & initialization \\
  & & 3_{\ }1_{0}5_{\ }2_{0}4_{\ } & right: 0 \rightarrow 0 \\
  & & 3_{1}1_{0}5_{\ }2_{0}4_{\ } & left: 0 \rightarrow 1 \\
  & & 3_{1}1_{0}5_{\ }2_{0}4_{1} & right: 1 \rightarrow 1 \\
  w(S) & = & 3_{1}1_{0}5_{2}2_{0}4_{1} & left: 1 \rightarrow 2 \\
\end{array}
\end{displaymath}
Filling the index in corresponding cell of the tableau $S$, We obtain
$i(S)$, the index tableau of $S$.
\begin{displaymath}
i(S) = {\def\lr#1{\multicolumn{1}{|@{\hspace{.6ex}}c@{\hspace{.6ex}}|}{\raisebox{-.3ex}{$#1$}}}
\raisebox{-.6ex}{$\begin{array}[]{*{3}c}\cline{1-2}
\lr{1}&\lr{2}\\\cline{1-3}
\lr{0}&\lr{0}&\lr{1}\\\cline{1-3}
\end{array}$}
}
\end{displaymath}
Now, using the tableaux $T$ and $i(S)$, cells of $T$ giving variable
index and corresponding cell of $i(S)$ giving exponent, we build
monomials $\x_T^{i(S)}$ as follow.
\begin{displaymath}
\begin{array}{c}
T = {\def\lr#1{\multicolumn{1}{|@{\hspace{.6ex}}c@{\hspace{.6ex}}|}{\raisebox{-.3ex}{$#1$}}}
\raisebox{-.6ex}{$\begin{array}[]{*{3}c}\cline{1-2}
\lr{2}&\lr{4}\\\cline{1-3}
\lr{1}&\lr{3}&\lr{5}\\\cline{1-3}
\end{array}$}
}
\qquad \qquad \qquad \qquad \qquad \qquad
i(S) = {\def\lr#1{\multicolumn{1}{|@{\hspace{.6ex}}c@{\hspace{.6ex}}|}{\raisebox{-.3ex}{$#1$}}}
\raisebox{-.6ex}{$\begin{array}[]{*{3}c}\cline{1-2}
\lr{1}&\lr{2}\\\cline{1-3}
\lr{0}&\lr{0}&\lr{1}\\\cline{1-3}
\end{array}$}
} \\
\x_T^{i(S)} = x_1^0 x_2^1 x_3^0 x_4^2 x_5^1
\end{array}
\end{displaymath}

Here are all the monomials in three variables.

\begin{displaymath}
\begin{array}{cccc}
S & T & i(S) & \x_T^{i(S)} \vspace{0.2cm} \\  
{\def\lr#1{\multicolumn{1}{|@{\hspace{.6ex}}c@{\hspace{.6ex}}|}{\raisebox{-.3ex}{$#1$}}}
\raisebox{-.6ex}{$\begin{array}[]{*{3}c}\cline{1-3}
\lr{1}&\lr{2}&\lr{3}\\\cline{1-3}
\end{array}$}}
&
{\def\lr#1{\multicolumn{1}{|@{\hspace{.6ex}}c@{\hspace{.6ex}}|}{\raisebox{-.3ex}{$#1$}}}
\raisebox{-.6ex}{$\begin{array}[]{*{3}c}\cline{1-3}
\lr{1}&\lr{2}&\lr{3}\\\cline{1-3}
\end{array}$}
}
& 
{\def\lr#1{\multicolumn{1}{|@{\hspace{.6ex}}c@{\hspace{.6ex}}|}{\raisebox{-.3ex}{$#1$}}}
\raisebox{-.6ex}{$\begin{array}[]{*{3}c}\cline{1-3}
\lr{0}&\lr{0}&\lr{0}\\\cline{1-3}
\end{array}$}}
& 1 \vspace{0.2cm} \\
{\def\lr#1{\multicolumn{1}{|@{\hspace{.6ex}}c@{\hspace{.6ex}}|}{\raisebox{-.3ex}{$#1$}}}
\raisebox{-.6ex}{$\begin{array}[]{*{2}c}\cline{1-1}
\lr{3}\\\cline{1-2}
\lr{1}&\lr{2}\\\cline{1-2}
\end{array}$}}
&
{\def\lr#1{\multicolumn{1}{|@{\hspace{.6ex}}c@{\hspace{.6ex}}|}{\raisebox{-.3ex}{$#1$}}}
\raisebox{-.6ex}{$\begin{array}[]{*{2}c}\cline{1-1}
\lr{3}\\\cline{1-2}
\lr{1}&\lr{2}\\\cline{1-2}
\end{array}$}}
&
{\def\lr#1{\multicolumn{1}{|@{\hspace{.6ex}}c@{\hspace{.6ex}}|}{\raisebox{-.3ex}{$#1$}}}
\raisebox{-.6ex}{$\begin{array}[]{*{2}c}\cline{1-1}
\lr{1}\\\cline{1-2}
\lr{0}&\lr{0}\\\cline{1-2}
\end{array}$}}
& x_3 \vspace{0.2cm} \\
{\def\lr#1{\multicolumn{1}{|@{\hspace{.6ex}}c@{\hspace{.6ex}}|}{\raisebox{-.3ex}{$#1$}}}
\raisebox{-.6ex}{$\begin{array}[]{*{2}c}\cline{1-1}
\lr{3}\\\cline{1-2}
\lr{1}&\lr{2}\\\cline{1-2}
\end{array}$}}
&
{\def\lr#1{\multicolumn{1}{|@{\hspace{.6ex}}c@{\hspace{.6ex}}|}{\raisebox{-.3ex}{$#1$}}}
\raisebox{-.6ex}{$\begin{array}[]{*{2}c}\cline{1-1}
\lr{2}\\\cline{1-2}
\lr{1}&\lr{3}\\\cline{1-2}
\end{array}$}}
&
{\def\lr#1{\multicolumn{1}{|@{\hspace{.6ex}}c@{\hspace{.6ex}}|}{\raisebox{-.3ex}{$#1$}}}
\raisebox{-.6ex}{$\begin{array}[]{*{2}c}\cline{1-1}
\lr{1}\\\cline{1-2}
\lr{0}&\lr{0}\\\cline{1-2}
\end{array}$}} 
& x_2 \vspace{0.2cm}  \\
{\def\lr#1{\multicolumn{1}{|@{\hspace{.6ex}}c@{\hspace{.6ex}}|}{\raisebox{-.3ex}{$#1$}}}
\raisebox{-.6ex}{$\begin{array}[]{*{2}c}\cline{1-1}
\lr{2}\\\cline{1-2}
\lr{1}&\lr{3}\\\cline{1-2}
\end{array}$}}
&
{\def\lr#1{\multicolumn{1}{|@{\hspace{.6ex}}c@{\hspace{.6ex}}|}{\raisebox{-.3ex}{$#1$}}}
\raisebox{-.6ex}{$\begin{array}[]{*{2}c}\cline{1-1}
\lr{3}\\\cline{1-2}
\lr{1}&\lr{2}\\\cline{1-2}
\end{array}$}} 
&
{\def\lr#1{\multicolumn{1}{|@{\hspace{.6ex}}c@{\hspace{.6ex}}|}{\raisebox{-.3ex}{$#1$}}}
\raisebox{-.6ex}{$\begin{array}[]{*{2}c}\cline{1-1}
\lr{1}\\\cline{1-2}
\lr{0}&\lr{1}\\\cline{1-2}
\end{array}$}} 
& x_2 x_3 \vspace{0.2cm} \\
{\def\lr#1{\multicolumn{1}{|@{\hspace{.6ex}}c@{\hspace{.6ex}}|}{\raisebox{-.3ex}{$#1$}}}
\raisebox{-.6ex}{$\begin{array}[]{*{2}c}\cline{1-1}
\lr{2}\\\cline{1-2}
\lr{1}&\lr{3}\\\cline{1-2}
\end{array}$}} 
&
{\def\lr#1{\multicolumn{1}{|@{\hspace{.6ex}}c@{\hspace{.6ex}}|}{\raisebox{-.3ex}{$#1$}}}
\raisebox{-.6ex}{$\begin{array}[]{*{2}c}\cline{1-1}
\lr{2}\\\cline{1-2}
\lr{1}&\lr{3}\\\cline{1-2}
\end{array}$}
}
& 
{\def\lr#1{\multicolumn{1}{|@{\hspace{.6ex}}c@{\hspace{.6ex}}|}{\raisebox{-.3ex}{$#1$}}}
\raisebox{-.6ex}{$\begin{array}[]{*{2}c}\cline{1-1}
\lr{1}\\\cline{1-2}
\lr{0}&\lr{1}\\\cline{1-2}
\end{array}$}}
& x_2 x_3 \vspace{0.2cm} \\
{\def\lr#1{\multicolumn{1}{|@{\hspace{.6ex}}c@{\hspace{.6ex}}|}{\raisebox{-.3ex}{$#1$}}}
\raisebox{-.6ex}{$\begin{array}[]{*{1}c}\cline{1-1}
\lr{3}\\\cline{1-1}
\lr{2}\\\cline{1-1}
\lr{1}\\\cline{1-1}
\end{array}$}}
&
{\def\lr#1{\multicolumn{1}{|@{\hspace{.6ex}}c@{\hspace{.6ex}}|}{\raisebox{-.3ex}{$#1$}}}
\raisebox{-.6ex}{$\begin{array}[]{*{1}c}\cline{1-1}
\lr{3}\\\cline{1-1}
\lr{2}\\\cline{1-1}
\lr{1}\\\cline{1-1}
\end{array}$}}
& 
{\def\lr#1{\multicolumn{1}{|@{\hspace{.6ex}}c@{\hspace{.6ex}}|}{\raisebox{-.3ex}{$#1$}}}
\raisebox{-.6ex}{$\begin{array}[]{*{1}c}\cline{1-1}
\lr{2}\\\cline{1-1}
\lr{1}\\\cline{1-1}
\lr{0}\\\cline{1-1}
\end{array}$}}
& x_2 x_3^2 \vspace{0.2cm} \\
\end{array}
\end{displaymath}

For $T$ a standard tableaux of shape $\lambda$, let $R(T)$ and $C(T)$
denote the row stabilizer and the column stabilizer of $T$
respectively and consider the Young symmetrizer
\begin{equation}
  \epsilon_{T} := \sum_{\sigma \in R(T)} \sum_{\tau \in C(T)} sign(\tau) \tau \sigma
\end{equation}
which is an element of the group algebra $\QQ[\sg[n]]$. We know define
the polynomial $F^S_T$ by
\begin{equation}
F^S_T(x_1, \dots , x_n) := \epsilon_T(\x_T^{i(S)}).
\end{equation}

\begin{theorem}
  Let $n$ a positive integer, the family of $n!$ polynomials $F^S_T$ for
  $S,T$ running over standard tableaux of the same shape form a basis
  of the $Sym(\x)$-module $\KK[\x]$.
\end{theorem}

Terasoma and Yamada proved it using the usual bilinear form in its
context : the divided difference associated to the longest element of
the symmetric group~\cite{MR1210951}.

In Three variable, here is the basis of $\KK[x_1, x_2, x_3]$ as a
$Sym(x_1, x_2, x_3)$-module.

\begin{displaymath}
\begin{array}{cccc}
S & T & \x_T^{i(S)} & F^S_T \vspace{0.2cm} \\  
{\def\lr#1{\multicolumn{1}{|@{\hspace{.6ex}}c@{\hspace{.6ex}}|}{\raisebox{-.3ex}{$#1$}}}
\raisebox{-.6ex}{$\begin{array}[]{*{3}c}\cline{1-3}
\lr{1}&\lr{2}&\lr{3}\\\cline{1-3}
\end{array}$}}
&
{\def\lr#1{\multicolumn{1}{|@{\hspace{.6ex}}c@{\hspace{.6ex}}|}{\raisebox{-.3ex}{$#1$}}}
\raisebox{-.6ex}{$\begin{array}[]{*{3}c}\cline{1-3}
\lr{1}&\lr{2}&\lr{3}\\\cline{1-3}
\end{array}$}
}
& 
1
& 6 \vspace{0.2cm} \\
{\def\lr#1{\multicolumn{1}{|@{\hspace{.6ex}}c@{\hspace{.6ex}}|}{\raisebox{-.3ex}{$#1$}}}
\raisebox{-.6ex}{$\begin{array}[]{*{2}c}\cline{1-1}
\lr{3}\\\cline{1-2}
\lr{1}&\lr{2}\\\cline{1-2}
\end{array}$}}
&
{\def\lr#1{\multicolumn{1}{|@{\hspace{.6ex}}c@{\hspace{.6ex}}|}{\raisebox{-.3ex}{$#1$}}}
\raisebox{-.6ex}{$\begin{array}[]{*{2}c}\cline{1-1}
\lr{3}\\\cline{1-2}
\lr{1}&\lr{2}\\\cline{1-2}
\end{array}$}}
&
x_3
& 2(x_3 - x_1) \vspace{0.2cm} \\
{\def\lr#1{\multicolumn{1}{|@{\hspace{.6ex}}c@{\hspace{.6ex}}|}{\raisebox{-.3ex}{$#1$}}}
\raisebox{-.6ex}{$\begin{array}[]{*{2}c}\cline{1-1}
\lr{3}\\\cline{1-2}
\lr{1}&\lr{2}\\\cline{1-2}
\end{array}$}}
&
{\def\lr#1{\multicolumn{1}{|@{\hspace{.6ex}}c@{\hspace{.6ex}}|}{\raisebox{-.3ex}{$#1$}}}
\raisebox{-.6ex}{$\begin{array}[]{*{2}c}\cline{1-1}
\lr{2}\\\cline{1-2}
\lr{1}&\lr{3}\\\cline{1-2}
\end{array}$}}
&
x_2
& 2(x_2 - x_1) \vspace{0.2cm}  \\
{\def\lr#1{\multicolumn{1}{|@{\hspace{.6ex}}c@{\hspace{.6ex}}|}{\raisebox{-.3ex}{$#1$}}}
\raisebox{-.6ex}{$\begin{array}[]{*{2}c}\cline{1-1}
\lr{2}\\\cline{1-2}
\lr{1}&\lr{3}\\\cline{1-2}
\end{array}$}}
&
{\def\lr#1{\multicolumn{1}{|@{\hspace{.6ex}}c@{\hspace{.6ex}}|}{\raisebox{-.3ex}{$#1$}}}
\raisebox{-.6ex}{$\begin{array}[]{*{2}c}\cline{1-1}
\lr{3}\\\cline{1-2}
\lr{1}&\lr{2}\\\cline{1-2}
\end{array}$}} 
&
x_2 x_3
& x_2 (x_3 - x_1) \vspace{0.2cm} \\
{\def\lr#1{\multicolumn{1}{|@{\hspace{.6ex}}c@{\hspace{.6ex}}|}{\raisebox{-.3ex}{$#1$}}}
\raisebox{-.6ex}{$\begin{array}[]{*{2}c}\cline{1-1}
\lr{2}\\\cline{1-2}
\lr{1}&\lr{3}\\\cline{1-2}
\end{array}$}} 
&
{\def\lr#1{\multicolumn{1}{|@{\hspace{.6ex}}c@{\hspace{.6ex}}|}{\raisebox{-.3ex}{$#1$}}}
\raisebox{-.6ex}{$\begin{array}[]{*{2}c}\cline{1-1}
\lr{2}\\\cline{1-2}
\lr{1}&\lr{3}\\\cline{1-2}
\end{array}$}
}
& 
x_2 x_3
& x_3 (x_2 - x_1) \vspace{0.2cm} \\
{\def\lr#1{\multicolumn{1}{|@{\hspace{.6ex}}c@{\hspace{.6ex}}|}{\raisebox{-.3ex}{$#1$}}}
\raisebox{-.6ex}{$\begin{array}[]{*{1}c}\cline{1-1}
\lr{3}\\\cline{1-1}
\lr{2}\\\cline{1-1}
\lr{1}\\\cline{1-1}
\end{array}$}}
&
{\def\lr#1{\multicolumn{1}{|@{\hspace{.6ex}}c@{\hspace{.6ex}}|}{\raisebox{-.3ex}{$#1$}}}
\raisebox{-.6ex}{$\begin{array}[]{*{1}c}\cline{1-1}
\lr{3}\\\cline{1-1}
\lr{2}\\\cline{1-1}
\lr{1}\\\cline{1-1}
\end{array}$}}
& 
x_2 x_3^2
& (x_3 - x_1)(x_3 - x_2)(x_2 - x_1) \vspace{0.2cm} \\
\end{array}
\end{displaymath}

We will know try to solve Problem~\ref{main-pb} by searching linear
combination of higher Specht polynomials stabilized by the action of a
permutation group.

\section{Combinatorial description of the invariant ring}

We now try to slice the invariant ring finer than degree by degree. As
irreducible representations of the symmetric group are homogeneous, we
will build format series mixing degree statistic and partitions.

\subsection{A refinement of the Moliens series}

Let $G \subset \sg[n]$ a permutation group. Any module $\sg[n]$-stable
is also $G$-stable, thus any representation of $\sg[n]$ is also a
representation of $G$. Usually, An irreducible representation of
$\sg[n]$ will not stay irreducible when restricted to $G$. We are
searching trivial representations of $G$ inside irreducible
representation of $\sg[n]$ and that can be done with a scalar product
of character.

For $G$ a permutation group, we will denote $\mathcal{C}(G)$ the set
of conjugacy classes of $G$. The usual scalar product between two
characters $\chi$ and $\psi$  of $G$ is given by
\begin{equation}
  \langle \chi , \psi \rangle  = \frac{1}{|G|} \sum_{C \in \mathcal{C}(G)}
  |C| \chi(\sigma) \psi(\sigma) \qquad (\sigma \text{ chosen arbitrary } \in C)
\end{equation}

\begin{proposition}
  Let $\lambda \vdash n$ a partition of the positive integer $n$. Let
  $G \subset \sg[n]$ a permutation group. The multiplicity of the
  trivial representation of $G$ inside the irreducible representation
  of $\sg[n]$ indexed by $\lambda$ is given by $m_{\lambda}(G,
  \sg[n])$ with 
\begin{equation}\label{trivial_multi}
  m_{\lambda}(G, \sg[n]) := \frac{1}{|G|} \sum_{C \in \mathcal{C}(G)}
  |C| M_{(\lambda, \text{cycle type}(\sigma))}, \qquad (\sigma \text{ chosen arbitrary } \in C)
\end{equation}
where $M_{(\lambda, cycle type(\sigma))}$ is the coefficient of the
character table of $\sg[n]$ indexed by the partitions $\lambda$ and
$\text{cycle type}(\sigma)$.

\begin{proof}
This just consist in using the usual scalar product of characters for
$G$ with the trivial character of $G$. Thus we can remark that value
of characters can be read on the character table of $\sg[n]$ because
conjugacy classes of $G$ are subset of conjugacy classes of
$\sg[n]$ and traces of matrices do not change when a representation of
$\sg[n]$ is viewed as a representation of $G$.
\end{proof}
\end{proposition}

\begin{definition}
  Let $G \subset \sg[n]$ a permutation group. Using a formal set of
  variable $\mathbf{t} = (t_{\lambda})_{\lambda \vdash n}$ indexed by
  partitions of $n$, we define the \emph{trivial multiplicities
    enumerator} $P(G, \mathbf{t})$ as follow
  \begin{equation}
    P(G, \mathbf{t}) := \sum_{\lambda \vdash n} m_{\lambda}(G, \sg[n]) t_{\lambda}
  \end{equation}
\end{definition}

$P(G, \mathbf{t})$ count the multiplicities of the trivial
representation of $G$ inside the irreducible representations of the
symmetric group of degree $n$ themselves indexed by the partitions of
the integer $n$.

For $G = \langle (1,2)(3,4), (1,4)(2,3) \rangle$.

\sageex{
\sagepromt G = PermutationGroup([[(1,2),(3,4)],[(1,4),(2,3)]]); G

\sageret{Permutation Group with generators [(1,2)(3,4), (1,4)(2,3)]}

\sagepromt for C in G.conjugacy\_classes():

\sagedots \quad card =  C.cardinality()

\sagedots \quad ct = Permutation(C.representative()).cycle\_type()
    
\sagedots \quad print card, ct

\sageret{1 [1, 1, 1, 1] \\
1 [2, 2] \\
1 [2, 2] \\
1 [2, 2] \\
}

\sagepromt trivial\_representations\_in\_symmetric\_representations(G)

\sageret{ \{[1, 1, 1, 1]: 1, [2, 1, 1]: 0, [2, 2]: 2, [3, 1]: 0, [4]: 1\} }
}
For this group, we thus have
\begin{displaymath}
P(G, \mathbf{t}) = t_{[1,1,1,1]} + 2t_{[2,2]} + t_{[4]}
\end{displaymath}

\begin{definition}
  Let $\lambda \vdash n$ a partition of a positive integer $n$ and $z$ a formal
  variable. The will denote $\phi(\lambda, z)$ the
  \emph{representation appearance polynomial} defined as follow
  \begin{equation}
    \phi(\lambda, z) := \sum_{T \in STab(\lambda)} z^{cocharge(T)}
  \end{equation}
  Where the sum run over all standard tableaux $T$ of shape $\lambda$. 
\end{definition}

$\phi(\lambda, z)$ make the link between the degree $z$ and the
irreducible representations of $\sg[n]$ isomorphic the abstract one
indexed by $\lambda$ appearing inside the quotient
$\KK[\x]/Sym^+(\x)$. $\phi(\lambda, 1) = |STab(\lambda)|$ give the
multiplicities of the irreducible representation indexed by
$\lambda$. Generally, if $\phi(\lambda, z)$ has coefficient an integer
$k$ for a term in $z^d$, this means that $k$ $Sym(\x)$-module
isomorphic to the irreducible representation of $\sg[n]$ indexed by
$\lambda$ can be built inside the graded quotient $\KK[\x]/Sym^+(\x)$
at degree $d$. The higher Specht Polynomials realize explicitly these
representations because the \emph{cocharge} is exactly the sum of the
entries of tableau $i(S)$ (or the degree of the corresponding Specht).

\begin{proposition}
  Let $G \subset \sg[n]$ a permutation group. The trivial
  multiplicities enumerator $P(G, t)$ and the Hilbert series $H(G, z)$ are related by
  \begin{equation}
    H(G, z) = \frac{P(G, t_{\lambda} \rightarrow \phi(\lambda, z))}{(1-z)(1-z^2) \cdots (1-z^n)}
  \end{equation}
\begin{proof}
  This result is a consequence of some statements about the
  combinatorics of standard tableaux. As discussed previously, once we
  have the dimension of $G$-trivial space inside each irreducible
  representation, it remains to know at which degree copies of
  irreducible representation lies into the quotient
  $\KK[\x]/Sym^+(\x)$. The \emph{cocharge} of standard tableaux is the
  right statistic partitioning the occurrences of $\sg[n]$-spaces
  along the degree.
\end{proof}
\end{proposition}

Back with the example $G = \langle (1,2)(3,4), (1,4)(2,3) \rangle$,

\sageex{
\sagepromt G = PermutationGroup([[(1,2),(3,4)],[(1,4),(2,3)]]); G

\sageret{Permutation Group with generators [(1,2)(3,4), (1,4)(2,3)]}

\sagepromt G.molien\_series()

\sageret{$(x^2 - x + 1)/(x^6 - 2x^5 - x^4 + 4x^3 - x^2 - 2x + 1)$}

\sagepromt S4 = SymmetricGroup(4); S4

\sageret{Symmetric group of order 4! as a permutation group}

\sagepromt S4.molien\_series()

\sageret{$1/(x^{10} - x^9 - x^8 + 2x^5 - x^2 - x + 1)$}

\sagepromt G.molien\_series() / S4.molien\_series()

\sageret{$x^6 + 2x^4 + 2x^2 + 1$}
}
On the other side, we had
\begin{equation}
P(G, \mathbf{t}) = t_{[1,1,1,1]} + 2t_{[2,2]} + t_{[4]}
\end{equation}

Let us now list standard tableaux of shape $[1,1,1,1], [2,2] and [4]$.

\begin{displaymath}
\begin{array}{ccc}
S & i(S) & cocharge(S) \vspace{0.2cm} \\
{\def\lr#1{\multicolumn{1}{|@{\hspace{.6ex}}c@{\hspace{.6ex}}|}{\raisebox{-.3ex}{$#1$}}}
\raisebox{-.6ex}{$\begin{array}[]{*{4}c}\cline{1-4}
\lr{1}&\lr{2}&\lr{3}&\lr{4}\\\cline{1-4}
\end{array}$}}
&
{\def\lr#1{\multicolumn{1}{|@{\hspace{.6ex}}c@{\hspace{.6ex}}|}{\raisebox{-.3ex}{$#1$}}}
\raisebox{-.6ex}{$\begin{array}[]{*{4}c}\cline{1-4}
\lr{0}&\lr{0}&\lr{0}&\lr{0}\\\cline{1-4}
\end{array}$}}
& 0 \vspace{0.2cm} \\
{\def\lr#1{\multicolumn{1}{|@{\hspace{.6ex}}c@{\hspace{.6ex}}|}{\raisebox{-.3ex}{$#1$}}}
\raisebox{-.6ex}{$\begin{array}[]{*{2}c}\cline{1-2}
\lr{3}&\lr{4}\\\cline{1-2}
\lr{1}&\lr{2}\\\cline{1-2}
\end{array}$}}
&
{\def\lr#1{\multicolumn{1}{|@{\hspace{.6ex}}c@{\hspace{.6ex}}|}{\raisebox{-.3ex}{$#1$}}}
\raisebox{-.6ex}{$\begin{array}[]{*{2}c}\cline{1-2}
\lr{1}&\lr{1}\\\cline{1-2}
\lr{0}&\lr{0}\\\cline{1-2}
\end{array}$}}
& 2 \vspace{0.2cm} \\
{\def\lr#1{\multicolumn{1}{|@{\hspace{.6ex}}c@{\hspace{.6ex}}|}{\raisebox{-.3ex}{$#1$}}}
\raisebox{-.6ex}{$\begin{array}[]{*{2}c}\cline{1-2}
\lr{2}&\lr{4}\\\cline{1-2}
\lr{1}&\lr{3}\\\cline{1-2}
\end{array}$}}
&
{\def\lr#1{\multicolumn{1}{|@{\hspace{.6ex}}c@{\hspace{.6ex}}|}{\raisebox{-.3ex}{$#1$}}}
\raisebox{-.6ex}{$\begin{array}[]{*{2}c}\cline{1-2}
\lr{1}&\lr{2}\\\cline{1-2}
\lr{0}&\lr{1}\\\cline{1-2}
\end{array}$}}
& 4 \vspace{0.2cm} \\
{\def\lr#1{\multicolumn{1}{|@{\hspace{.6ex}}c@{\hspace{.6ex}}|}{\raisebox{-.3ex}{$#1$}}}
\raisebox{-.6ex}{$\begin{array}[]{*{1}c}\cline{1-1}
\lr{4}\\\cline{1-1}
\lr{3}\\\cline{1-1}
\lr{2}\\\cline{1-1}
\lr{1}\\\cline{1-1}
\end{array}$}}
&
{\def\lr#1{\multicolumn{1}{|@{\hspace{.6ex}}c@{\hspace{.6ex}}|}{\raisebox{-.3ex}{$#1$}}}
\raisebox{-.6ex}{$\begin{array}[]{*{1}c}\cline{1-1}
\lr{3}\\\cline{1-1}
\lr{2}\\\cline{1-1}
\lr{1}\\\cline{1-1}
\lr{0}\\\cline{1-1}
\end{array}$}}
& 6 \vspace{0.2cm} \\
\end{array}
\end{displaymath}
We have $\phi([1,1,1,1], z) = 1$, $\phi([2,2], z) = z^2 + z^4$ and
$\phi([4], z) = z^6$. Injecting these evaluations inside the trivial
multiplicities enumerator, we recover
\begin{equation}
P(G, t_{\lambda} \rightarrow \phi(\lambda, z)) = 1 + 2z^2 + 2z^4 + z^6
\end{equation}
which give the number of secondary invariants degree by degree when
elementary symmetric polynomials are taken as primary invariants. This
polynomials is also the quotient of the two Hilbert series that can be
computed in GAP~\cite{GAP} using Molien's Formula.

\subsection{Secondary invariants built from higher Specht polynomial}

Let $G \in \sg[n]$ a permutation group, $\lambda \vdash n$ a partition
of $n$. Let us suppose that we have calculated $m_{\lambda}(G,
\sg[n])$ and we now want to build explicitly the secondary invariant
polynomials. We are in the case in which we have an homogeneous
$G$-stable space inside which we want to construct a finite and known
number of independent invariant polynomials under the action of $G$.

The usual way to dealt with this problem is to built an explicit
family spanning the concerned space by generating polynomials forming
a basis. Then basis element by basis element, we apply the Reynolds
operator and some linear algebra to get a free family of the wanted
dimension. knowing this dimension give a stopping criteria often very
important since computations are extremely heavy even for small number
of variables.

In our context, even the usual approach would work, as permutation are
often given by a list of generators, we can even forget the Reynolds
operator.

\begin{proposition}
Let $G \in \sg[n]$ a permutation group given by some generator : $G =
\langle \sigma_1 , \dots , \sigma_r \rangle$. Let $\lambda \vdash n$ a
partition. The $G$-trivial abstract space inside the abstract
representation of $\sg[n]$ indexed by $\lambda$ is given by the
intersection of the eigenspace of the representation matrices of
$\sigma_1 , \dots , \sigma_r$ associated the eigenvalue $1$.
\begin{proof}
We view here representation indexed by $\lambda$ as the formal free
module generated by the standard tableaux of shape $\lambda$. We know
that a subspace of dimension $m_{\lambda}(G, \sg[n])$ have its
elements invariant under the action of $G$. Being invariant under $G$
is equivalent to be stabilized by the Reynolds operator but is also
equivalent (by definition in fact) to be stabilized by the action of
the generators of $G$. Since we are working inside a representation of
$\sg[n]$, each permutation have an associated matrix and the kernel of
this matrix characterize the formal subspace stabilized

\end{proof}
\end{proposition}

\section{Algorithm building secondary invariants}

We now present an effective algorithm exploiting the approach using
the representation of the symmetric group.

Computation dependencies :
\begin{itemize}
  \item Character table of the symmetric group
  \item Conjugacy classes of the group $G$ (cardinalities and representatives)
  \item Matrices of irreducible representation of the symmetric group
  \item Some linear algebra
\end{itemize}

\begin{algorithm}
  \title{Compute secondary invariant using $\sg[n]$ representations}
  \label{sec_inv_algo}
  \textbf{Input :} $\sigma_1, \sigma_2, \dots \sigma_n$ a set of permutations of size $n$ generating a group $G$.

\begin{displaymath}
%  \begin{small}
  \begin{array}{l}
    \textbf{def }\text{SecondaryInvariants($\sigma_1, \sigma_2, \dots \sigma_n$)} : \\
    \ind G \leftarrow PermutationGroup(\sigma_1, \sigma_2, \dots \sigma_n) \\
    \ind \sg[n] \leftarrow SymmetricGroup(n) \\
    \ind S \leftarrow \{\} \\
    \ind \textbf{for } \lambda \in  Partition(n) : \\
    \ind \ind \textbf{if } m_{\lambda}(G, \sg[n]) \neq 0 : \\
    \ind \ind \ind G \leftarrow VectorSpace(\mathbb{Q}, f_{\lambda}) \\
    \ind \ind \ind \textbf{for } i \in \{1, 2, \dots , n\}: \\
    \ind \ind \ind \ind G \leftarrow G \cap kernel(M_{\lambda}(\sigma_i) - Id) \qquad \qquad \qquad \qquad \qquad \qquad \quad \\
    \ind \ind \ind abstract\_secondary \leftarrow basis(G) \\ 
    \ind \ind \ind \textbf{for } \mu \in StandardTableaux(\lambda) : \\
    \ind \ind \ind \ind \textbf{for } P \in abstract\_secondary : \\    
    \ind \ind \ind S \leftarrow S \cup HigherSpecthPolynomials(P, \mu) \\
    \ind \textbf{return } S \\
  \end{array}
%  \end{small}
\end{displaymath}

\end{algorithm}

The returned set is composed by linear combinations of higher Specht
polynomials. These polynomials can be easily evaluated but, as they
contains a lot of Vandermonde factors, there expansion on a set of $n$
formal variables is an heavy computation.

\subsection{A large trace of the algorithm}

A teen years computational challenge consist in computing a generating
family of the ring of invariants of the group acting on the edges of
graphs over $5$ nodes. This group is a subgroup of the symmetric group
of degree $10 = \binom{5}{2}$ and has for cardinality $5! = 120$. As
far as we know, no computer algebra system has already handle a such
computation. After $12$ hours of computation, Magma, Singular and the
evaluation approach written in Sage did not finish (in Sage, more
precisely, around 5 percent of this computation(linear algebra) was
done after $24$ hours). We tried our approach on this group and we got
the following verbose :

The trace have should be read has the following pattern
\begin{footnotesize}
\begin{verbatim}
[partition]  ambient dimension -->  (number of standard tableaux for this shape)
rank in S_n repr : (dimension of the G-trivial space)
\end{verbatim}
\end{footnotesize}

\begin{footnotesize}
\begin{verbatim}
sage: load("invariants.py")
sage: G = TransitiveGroup(10,12)
sage: Specht_basis_of_trivial_representations(G, verbose=True)
[3, 2, 2, 1, 1, 1]  ambient dimension -->  315
rank in S_n repr :  2
[6, 1, 1, 1, 1]  ambient dimension -->  126
rank in S_n repr :  3
[6, 4]  ambient dimension -->  90
rank in S_n repr :  3
[4, 3, 1, 1, 1]  ambient dimension -->  525
rank in S_n repr :  5
[5, 2, 2, 1]  ambient dimension -->  525
rank in S_n repr :  4
[5, 2, 1, 1, 1]  ambient dimension -->  448
rank in S_n repr :  4
[3, 2, 1, 1, 1, 1, 1]  ambient dimension -->  160
rank in S_n repr :  1
[6, 2, 1, 1]  ambient dimension -->  350
rank in S_n repr :  2
[4, 4, 2]  ambient dimension -->  252
rank in S_n repr :  5
[2, 1, 1, 1, 1, 1, 1, 1, 1]  ambient dimension -->  9
rank in S_n repr :  1
[3, 3, 2, 2]  ambient dimension -->  252
rank in S_n repr :  2
[4, 4, 1, 1]  ambient dimension -->  300
rank in S_n repr :  2
[4, 2, 2, 2]  ambient dimension -->  300
rank in S_n repr :  5
[4, 2, 2, 1, 1]  ambient dimension -->  567
rank in S_n repr :  3
[2, 2, 2, 1, 1, 1, 1]  ambient dimension -->  75
rank in S_n repr :  2
[3, 2, 2, 2, 1]  ambient dimension -->  288
rank in S_n repr :  3
[4, 2, 1, 1, 1, 1]  ambient dimension -->  350
rank in S_n repr :  3
[3, 1, 1, 1, 1, 1, 1, 1]  ambient dimension -->  36
rank in S_n repr :  1
[7, 1, 1, 1]  ambient dimension -->  84
rank in S_n repr :  1
[5, 1, 1, 1, 1, 1]  ambient dimension -->  126
rank in S_n repr :  3
[2, 2, 2, 2, 2]  ambient dimension -->  42
rank in S_n repr :  3
[6, 3, 1]  ambient dimension -->  315
rank in S_n repr :  1
[8, 2]  ambient dimension -->  35
rank in S_n repr :  2
[3, 3, 3, 1]  ambient dimension -->  210
rank in S_n repr :  2
[3, 3, 1, 1, 1, 1]  ambient dimension -->  225
rank in S_n repr :  1
[5, 4, 1]  ambient dimension -->  288
rank in S_n repr :  3
[5, 3, 2]  ambient dimension -->  450
rank in S_n repr :  3
[10]  ambient dimension -->  1
rank in S_n repr :  1
[4, 3, 2, 1]  ambient dimension -->  768
rank in S_n repr :  6
[7, 2, 1]  ambient dimension -->  160
rank in S_n repr :  1
[6, 2, 2]  ambient dimension -->  225
rank in S_n repr :  3
[5, 3, 1, 1]  ambient dimension -->  567
rank in S_n repr :  5
[3, 3, 2, 1, 1]  ambient dimension -->  450
rank in S_n repr :  4
total :  30240
n! / |G| :  30240
TOTAL CPU TIME :  414.837207
\end{verbatim}
\end{footnotesize}

Our algorithm took $414$ seconds to generated the $30240$ secondary
invariants as linear combinations of higher Specht polynomials. We
still believe that the computation of a couple primary and secondary
invariants for this group is unreachable for Magma, Singular and the
evaluation approach in less than $24$ hours.

For example, inside the symmetric representation associated with the
partition $(4,3,2,1)$. The algorithm making the Gauss reduction build a
space of dimension $6$ inside an ambient space of dimension
$768$. Using the Hook-length formula, on can check that there exist $768$ standard tableaux of shape $(4,3,2,1)$
\begin{equation}
{\def\lr#1{\multicolumn{1}{|@{\hspace{.6ex}}c@{\hspace{.6ex}}|}{\raisebox{-.3ex}{$#1$}}}
\raisebox{-.6ex}{$\begin{array}[]{*{4}c}\cline{1-1}
\lr{1}\\\cline{1-2}
\lr{3}&\lr{1}\\\cline{1-3}
\lr{5}&\lr{3}&\lr{1}\\\cline{1-4}
\lr{7}&\lr{5}&\lr{3}&\lr{1}\\\cline{1-4}
\end{array}$}
} \qquad \frac{10!}{7 \cdot 5 \cdot 5 \cdot 3 \cdot 3 \cdot 3} = 768
\end{equation}

\subsection{Implementation details}

Our approach involves a lot a technology around representation theory
and symbolic computation, most of the required prerequisites are
available in current computer algebras system. For our problem,
dependencies have been implemented in \gap~\cite{GAP} or
\Sage~\cite{Sage}. We implemented a first version in \Sage (400 lines
of code with tests and documentation) to test the efficienty of this
approach.

Here are the main steps of computation.

\begin{itemize}
  \item Compute the character table of the symmetric group. \\ 
    This task is done by \gap~\cite{GAP}. The Murnaghan-Nakayama rule
    is an example of working algorithm. This is a small part of
    computations and it can also be precomputed and store for all
    small symmetric group (for degree at most $20$).

  \item Enumerate cardinality and a representative of each conjugacy class of $G$. \\
    Also handled by \gap~\cite{GAP} and interfaced in
    \Sage~\cite{Sage}. The general problem of enumerating conjugacy
    classes of finite group is not simple but there are efficient
    algorithms for permutation group.

  \item Compute the trivial multiplicities enumerator of $G$. \\
    We computed a short function which just iterate over all partitions
    and compute a simple scalar product (see Formula~\ref{trivial_multi}). 

  \item Calculate the matrices of $\sigma_1, \sigma_2, \dots , \sigma_n$ inside the irreducible representations of $\sg[n]$. \\
    Once we have identified a $G$-trivial space inside an irreducible
    of $\sg[n]$, we build the matrices of the generators of $G$ using
    the \Sage~\cite{Sage} object :
    \texttt{SymmetricGroupRepresentation} which admits as argument a
    partition. The returned object is able to build matrices of the
    abstract representation from a permutation given in line notation. 

  \item Compute the intersection of stabilized subspace of these matrices. \\
    We use in a loop the method \texttt{intersection} of \Sage
    \texttt{VectorSpace} and the method \texttt{kernel} of \Sage
    \texttt{Matrix} class. Our linear algebra is thus handle by
    \Sage~\cite{Sage}.

  \item Transcript abstract combinations of standard tableaux in term of higher Specht polynomials. \\
    \Sage and especially \sagecombinat~\cite{Sage-Combinat} contains
    the combinatorial family of standard tableaux. A lot of
    combinatorial statistic are availlable on such objects like the
    cocharge. We implemented the higher Specht polynomials as \Sage
    polynomials computed from a pair of \Sage
    \texttt{StandardTableaux} of the same shape. 
\end{itemize}

Our code is completely not optimized. Since most of computations in
invariants (minimal generating set, secondary invariants, Hironaka
decomposition) present large theoretical exponential complexity, we
believe that, before code source refinement, we need new algorithms
with better asymptotically behavior. However, a huge factor of
execution time can easily be wined from our code with cache remember, better
interface between \gap and \Sage, parralelisation (each computer core
can echelonize matrices for a special irreducible representation), we
wanted to benchmarks our approach using the current available (and
open source) technology. Therefore, anyone (even a student) can
reproduce our experience rapidly and for free.

\subsection{Complexity}

The rich literature about effective invariant theory does not provide
a lot a fin complexity bound for algorithms. Gröbner basis admit very
general complexity bounds (in worst case $2^{2^{O(n)}}$ for n variables)
which appears to be overestimated compared to their effective
behavior. The thesis~\cite{Borie.2011.Thesis} of the author present an
evaluation approach to compute the invariants inside a quotient of a
reduced dimension. The algorithm computing secondary with this
technique has a complexity in $O((n!)^2 +\frac{(n!)^3}{|G|^2}
)$. Using the representation of the symmetric group, it is still very
hard to establish a fin bound. However, we can produce better bounds.

\begin{theorem}
Let $G$ a permutation group, subgroup of $\sg[n]$, given by $r$
generators. The complexity of the linear algebra reduction in
algorithm computing the secondary invariants of $G$ in
Section~\ref{sec_inv_algo} is bounded by
\begin{equation}
r \cdot ( \sum_{\lambda \vdash n} m_{\lambda}(|G|, \sg[n]) (f_{\lambda})^2 )
\end{equation}
where $f_{\lambda}$ is the number of standard tableaux of shape
$\lambda$ and $p(n)$ is the number of partition of size $n$.
\begin{proof}
This is just a straightforward counting of the reductions the $r$
matrices of permutation for each irreducible representation of the
symmetric group. For each partition $\lambda$ of $n$, we have to
construct a free family of $m_{\lambda}(|G|, \sg[n])$ vectors
(i.e. the rank) inside a space of dimension $f^{\lambda}$ (i.e. number
of indeterminates) from $r \cdot f^{\lambda}$ equations (the vertical
concatenation of the $r$ matrices). Summing these operations give the
announced bound.
\end{proof}
\end{theorem}

\begin{corollary}
Let $G$ a permutation group, subgroup of $\sg[n]$, given by $r$
generators. The complexity of the algorithm computing the secondary
invariants in Section~\ref{sec_inv_algo} has a complexity in $O(r \cdot (n!)^{\frac{3}{2}})$.
\begin{proof}
Let us denote $f_{max} := \displaystyle\max_{\lambda \vdash n} \{ f^{\lambda} \}$. We thus have
\begin{equation}
r \cdot ( \sum_{\lambda \vdash n} m_{\lambda}(|G|, \sg[n]) (f_{\lambda})^2 ) \leqslant r \cdot ( \sum_{\lambda \vdash n} f_{max} (f_{\lambda})^2 )
\end{equation}
since, at worse, the irreducible representation of $\sg[n]$ is
composed only by element stabilized point wise by $G$. The $f_{max}$
term can go out the sum and using formula~\ref{sum_f_lambda} we get
the bound
\begin{equation}
r \cdot f_{max} \cdot  n!
\end{equation}
Using the triangular inequality, we roughly have $f_{max}^2 \leqslant n!$ and that give the result.
\end{proof}
\end{corollary}

\section{Further developments}

Even the impressive algorithmic efficiently of this approach using the
representation of the symmetric group, it has also the advantage of
putting a lot of combinatrics inside this problem often classed inside
effective algebraic geometry.

\section*{acknowledgements}
\label{sec:ack}

This research was driven by computer exploration using the open-source
mathematical software \Sage~\cite{Sage}. In particular, we perused its
algebraic combinatorics features developed by the \sagecombinat
community~\cite{Sage-Combinat}, as well as its group theoretical
features provided by \gap~\cite{GAP}.

\bibliographystyle{abbrv}
\bibliography{main}

\end{document}